# ВАРИАНТ ПОСТРОЕНИЯ И ОБОСНОВАНИЯ ОБОБЩЕННОЙ ФОРМУЛЫ ИТО-ВЕНТЦЕЛЯ

## Дубко В.А.

**Аннотация**

В работе рассмотрен один из вариантов обоснования обобщенной формулы Ито-Вентцеля. Установлена взаимосвязь между различными представлениями обобщенной формулы Ито-Вентцеля.

Ключевые слова: случайная пуассоновская мера, процесс Винера, обобщенное уравнение Ито.

## Введение.

В работе [2] был предложен алгоритм построения, обобщающий формулу Ито-Вентцеля, на случай обобщенных уравнений Ито и опирающийся на уравнения для ядер интегральных инвариантов. В дальнейшем, появились статьи, использующие различные подходы при ее обосновании, привлекающие, в том числе, и различные представления обобщенных уравнений Ито [4].

Внимание к этой формуле вызвано тем, что она является одним из важных элементов формирования стохастического анализа и его применения. Так, например, в [3] эта формула стала основой для доказательства теорем об уравнениях для стохастических первых интегралов, стохастических ядер обобщенных уравнений Ито.

В данной работе, представим вариант получения и обоснования обобщенной формулы Ито-Вентцеля и установим связи между коэффициентами различных представлений этой формулы.

## Обобщенная формула Ито-Вентцеля.

Пусть $x(t) \in \mathbb{R}^n$ решение обобщенного уравнения Ито:

$$dx(t) = a(t)dt + b_k(t)dw_k(t) + \int g(t;\gamma)\nu(dt;d\gamma)$$

$$x(t) = x(t;z)|_{t=0} = z \in \mathbb{R}^n, \quad k = \overline{1,m} \tag{1}$$

где $w(t)$ — $m$-мерный винеровский процесс, $\nu(t;\Delta\gamma)$ — однородная по $t$ случайная мера Пуассона, $\int \sim \int_{R(\gamma)}$, и по индексам, встречающимся дважды, здесь и далее, ведется суммирование.

Относительно коэффициентов (1) полагаем, что они согласовываются с требованиями, обеспечивающими существование и единственность решений [1].

Теорема. Пусть $F(t;x), (t, x) \in [0, T] \times \mathbb{R}^n$ - скалярная функция подчиненная стохастическому уравнению:

$$dF(t;x) = Q(t;x)dt + D_k(t;x)dw_k(t) + \int G(t;x;\gamma)\nu(dt;d\gamma), \quad (2)$$

$$F(t;x)|_{t=0} = F(x) \in C_0^2$$

**L.a).** $Q(t;x), D_k(t;x), G(t;x;\gamma) \in \mathbb{R}$ неупреждающие, в общем, случайные функции, измеримые относительно того же потока $\sigma-$алгебр $\mathcal{F}_t$, что и $w(t)$, $\nu(t;A), \forall A \in \boldsymbol{B}$ – фиксированной борелевской $\sigma-$алгебры.

**L.в).** С вероятностью единица, $Q(t;x), D_k(t;x), G(t;x;\gamma)$, вместе со своими вторыми частными производными по компонентам $x$, непрерывны и ограничены $\forall t, x$.

**L.c).** $\int |\nabla^\beta G(t;x;\gamma)|^\alpha \prod(d\gamma) < \infty, \forall x; \alpha = 1,2; \nabla^\beta$ - условное обозначение частных производных по компонентам параметра $x$, порядка $\beta = 0,1,2$.

Тогда, если $x(t)$ случайный процесс, подчиненный системе уравнений (1), то справедлива формула:

$$dF(t;x(t)) = Q(t;x(t))dt + D_k(t;x(t))dw_k + b_{i,k}(t)\frac{\partial}{\partial x_i}F(t;x(t))dw_k +$$
$$+ [a_i(t)\frac{\partial}{\partial x_i}F(t;x(t)) + \frac{1}{2}b_{i,k}(t)b_{j,k}(t)\frac{\partial^2 F(t;x(t))}{\partial x_i \partial x_j} + \quad (3)$$
$$+ b_{i,k}(t)\frac{\partial}{\partial x_i}D_k(t;x(t))]dt +$$
$$+ \int[(F(t;x(t)+g(t;\gamma))-F(t;x(t)))]\nu(dt;d\gamma) + \int G(t;x(t)+g(t;\gamma);\gamma)\nu(dt;d\gamma).$$

Доказательство. Отметим, что при указанных ограничениях **L)** возможно дифференцирование $F(t;x)$ по компонентам параметра $x$. Соответствующие решения для производных и смешанных производных, построенных на основе (2), существуют, единственны и приводят к таким оценкам и свойствам:

$L_1$) $P(|\nabla^\beta F(t;x)| > \varepsilon^{-1}) \leq \varepsilon^2 C, C -$ некоторая постоянная,

$L_2$) $\underset{|x| \to 0}{\text{l.i.m.}} | F(t;x+y) - F(t;y)| = 0, \underset{|x| \to 0}{\text{l.i.m.}} | F'_{x_i}(t;x+y) - F'_{x_i}(t;y)| = 0,$

$\underset{|x| \to 0}{\text{l.i.m.}} | F''_{x_i x_j}(t;x+y) - F''_{x_i x_j}(t;y)| = 0, \forall (t, y) \in [0, T] \times \mathbb{R}^n.$

Заметим, что уравнение (2), в силу произвольности $x$, эквивалентно следующему:

$$d_t F(t;x+y) = Q(t;x+y)dt + D_k(t;x+y)dw_k(t) + \int G(t;x+y;\gamma)\nu(dt;d\gamma) \quad (4)$$

Пусть $f(x)$ детерминированная, непрерывная вместе со своими частными производными, вплоть до вторых включительно и

K). $\int_{\mathbb{R}^n} |f(x)| d\Gamma(x), \int_{\mathbb{R}^n} |\frac{\partial f(x)}{\partial x_i}| d\Gamma(x), \int_{\mathbb{R}^n} |\frac{\partial^2 f(x)}{\partial x_j \partial x_i}| d\Gamma(x) \leq const,$

$\lim_{|x|\to\infty} f(x) = 0; \lim_{|x|\to\infty} \frac{\partial f(x)}{\partial x_i} = 0.$

Воспользовавшись обобщенной формулой Ито, уравнениями (1) и (2), можем убедиться, что

$$dF(t;y)f(y-x(t;z)) = f(y-x(t;z))[Q(t;y)dt + D_k(t;y)dw_k(t)] - \\ -F(t;y)f_{y_i}(y-x(t;z))[a_i(t)dt + b_{i,k}(t)dw_k(t)] + \\ +0{,}5 \cdot b_{i,k}(t)b_{j,k}(t))f_{y_i,y_j}(y-x(t;z)) - D_k(t;y)b_{i,k}(t)f_{y_i}(y-x(t;z))dt + \\ +\int [f(y-x(t;z)-g(t;\gamma))\{F(t;y)+G(t;y;\gamma)\} - f(y-x(t;z))F(t;y)]\nu(dt;d\gamma) \quad (5)$$

С учетом требований к коэффициентам уравнения (2),(3), свойств L) и ограничений для $f(x)$, можно определить интегралы по $y$ $\forall x(t;z),t,\gamma$:

$\int_{\mathbb{R}^n} d\Gamma(y)F(t;y)f(y-x(t;z)), \int_{\mathbb{R}^n} d\Gamma(y)f(y-x(t;z))Q(t;y), (d\Gamma(y) = \prod_{i=1}^n dy_i),$

$\int_{\mathbb{R}^n} d\Gamma(y)f_{y_i}(y-x(t;z))D_k(t;y), \int_{\mathbb{R}^n} d\Gamma(y)f(y-x(t;z))D_k(t;y),$

$\int_{\mathbb{R}^n} d\Gamma(y)f_{y_i,y_j}(y-x(t;z)), \int_{\mathbb{R}^n} d\Gamma(y)F(t;y)f_{y_i}(y-x(t;z)),$

$\int_{\mathbb{R}^n} d\Gamma(y)[f(y-x(t;z)-g(t;\gamma))\{F(t;y)+G(t;y;\gamma)\} - f(y-x(t;z))F(t;y)]$,

трактуя их как среднеквадратичные пределы сумм Римана.

Проинтегрируем (5) по пространству параметра $y$, и, затем, выполнив интегрирование по частям, перенесем производные на $F(t;y)$. В результате приходим к равенству:

$$\int_{\mathbb{R}^n} d\Gamma(y)[F(t;y)f(y-x(t;z)) - F(y)f(y-z)] = \\ = \int_0^\tau dt\{\int_{\mathbb{R}^n} d\Gamma(y)f(y-x(t;z))\{[Q(t;y) + F_{y_i}(t;y)a_i(t) + D_{y_i k}(t;y)b_{i,k}(t) + \\ +\frac{1}{2} \cdot F_{y_i,y_j}(t;y)b_{i,k}(t)b_{j,k}(t) + D_{y_i k}(t;y)b_{i,k}(t)]dt + [D_k(t;y) + b_{i,k}(t)]dw_k(t)\} + \\ + \int_0^\tau \int \nu(dt;d\gamma)\int_{\mathbb{R}^n} d\Gamma(y)\{[f(y-x(t;z)-g(t;\gamma))F(t;y) - f(y-x(t;z))F(t;y)] + \\ + f(y-x(t;z)-g(t;\gamma))G(t;y;\gamma)\}. \quad (6)$$

Произведем в подынтегральных выражениях замену переменных:
а) $y - x(t;z) = \tilde{y}$, в интегралах связанных $t$ и $w_k(t)$, и переобозначим $\tilde{y} \to y$;
в) $y - x(t;z) - g(t;\gamma) = \tilde{y}$, в интегралах, связанных с пуассоновской мерой, и переобозначим $\tilde{y} \to y$.

В результате приходим к равенству:
$$\int_{\mathbb{R}^n} d\Gamma(y)f(y)[F(\tau;y+x(\tau;z))-F(y+z))] =$$
$$= \int_{\mathbb{R}^n} d\Gamma(y)f(y)\{\int_0^\tau ([Q(t;y+x(t;z))dt + D_k(t;y+x(t;z))dw_k(t)] +$$
$$+ F_{y_i}(t;y+x(t;z))[a_i(t)dt + b_{i,k}(t)dw_k(t)] +$$
$$+ 0{,}5 \cdot F_{y_i,y_j}(t;y+x(t;z))b_{i,k}(t)b_{j,k}(t))dt + D_{y_i k}(t;y+x(t;z))b_{i,k}(t)dt +$$
$$+ \int [(F(t;x(t;z)+g(t;\gamma)+y) - F(t;x(t;z)+y)]\nu(dt;d\gamma) +$$
$$+ \int G(t;x(t;z)+g(t;\gamma)+y;\gamma)\nu(dt;d\gamma))\}.$$

Перейдем, затем, от производных по $y_i, y_j$ к производным по $x_i, x_j$:
$$\int_{\mathbb{R}^n} d\Gamma(y)f(y)([F(\tau;y+x(\tau;z))-F(y+z))] -$$
$$- \int_0^\tau \{[Q(t;y+x(t;z))dt + D_k(t;y+x(t;z))dw_k(t)] +$$
$$+ F_{x_i}(t;y+x(t;z))[a_i(t)dt + b_{i,k}(t))dw_k(t)] + D_{x_i k}(t;y+x(t;z))b_{i,k}(t))dt +$$
$$+ \frac{1}{2}\cdot b_{i,k}(t)b_{j,k}(t)F_{x_i,x_j}(t;y+x(t;z))dt + b_{i,k}(t)D_{x_i k}(t;y+x(t;z))dt +$$
$$+ \int[(F(t;x(t;z)+y+g(t;\gamma)) - F(t;x(t;z)+y]\nu(dt;d\gamma) +$$
$$+ \int G(t;x(t;z)+g(t;\gamma)+y;\gamma)\nu(dt;d\gamma)\}) = 0.$$

В силу произвольности $f(y)$ из класса K), следствий L), гарантирующих существование и единственность, на каждой из траекторий $x(t;z)$, $\forall y$, функций и коэффициентов, зависящих от $(x(t;z)+y)$ в равенстве, приходим к выводу, что оно будет выполнено, только тогда, когда $\forall y, z$:
$$dF(t;x(t;z)+y) = Q(t;x(t;z)+y)dt + D_k(t;x(t;z)+y)dw_k +$$
$$+ b_{i,k}(t)\frac{\partial}{\partial x_i}F(t;x(t;z)+y)dw_k + [a_i(t)\frac{\partial}{\partial x_i}F(t;x(t;z)+y) +$$
$$+ \frac{1}{2}b_{i,k}(t)b_{j,k}(t)\frac{\partial^2 F(t;x(t;z)+y)}{\partial x_i \partial x_j} + b_{i,k}(t)\frac{\partial}{\partial x_i}D_k(t;x(t;z)+y)]dt +$$
$$+ \int G(t;x(t)+g(t;\gamma)+y;\gamma)\nu(dt;d\gamma) +$$
$$+ \int[(F(t;x(t)+g(t;\gamma)+y) - F(t;x(t;z))+y]\nu(dt;d\gamma).$$

С учетом замечания (4), выбирая $y = 0$, и приходим к формуле (3), одному из представлений обобщенной формулы Ито-Вентцеля.

## О различных представлениях обобщенной формулы Ито-Вентцеля.

Вид уравнений (3), зависят от вида уравнений для $x(t;z)$ и (2). Остановимся на некоторых различных вариантах представления формулы Ито-Вентцеля [3] и установим взаимосвязь между ними.

Рассмотрим два представления обобщенного уравнения Ито:

$$dx(t) = a^1(t)dt + b_k(t)dw_k(t) + \int g(t;\gamma)v(dt;d\gamma), \qquad (7)$$

$$dx(t) = \tilde{a}(t)dt + b_k(t)dw_k(t) + \int g(t;\gamma)\tilde{v}(dt;d\gamma). \qquad (8)$$

где $\tilde{v}(dt;d\gamma) = v(dt;d\gamma) - dt\prod(d\gamma)$ - центрированная пуассоновская мера,

$$\int_0^T dt \int |g(t;\gamma)|\prod(d\gamma) < \infty, \quad \mathrm{M}[v(dt;d\gamma)] = dt\prod(d\gamma).$$

Пусть $x(t)$ – решение системы (8), а $F(t;x)$ - системы (2). Перейдем от (8), к представлению (7). Это приведет к тому, что коэффициенты $a^1(t)$ в (7) должны быть представлены в виде:

$$a^1(t) = \tilde{a}(t) - \int g(t;\gamma)\prod(d\gamma).$$

Пусть, теперь, $x(t)$ – решение системы (1), а $F(t;x)$ - системы

$$dF(t;x) = \tilde{Q}(t;x)dt + D_k(t;x)dw_k(t) + \int G(t;x;\gamma)\tilde{v}(dt;d\gamma) \qquad (9)$$

Выполним переход от (9) к (2):

$$dF(t;x) = [\tilde{Q}(t;x) - \int G(t;x;\gamma)\prod(d\gamma)]dt + D_k(t;x)dw_k(t) + \int G(t;x;\gamma)v(dt;d\gamma)$$

Для того же, чтобы получить правило нахождения $dF(t;x(t))$, обусловленного системой (1), (9), в (2) необходимо трактовать $Q(t;x)$ как выражение:

$$Q(t;x) = \tilde{Q}(t;x) - \int G(t;x;\gamma)\prod(d\gamma)].$$

Если рассматривать как исходные системы (8), (9), то в (1), (2) требуется перейти к коэффициентам:

$$a(t) = \tilde{a}(t) - \int g_i(t;\gamma)\prod(d\gamma), \quad Q(t;x) = \tilde{Q}(t;x) - \int G(t;x;\gamma)\prod(d\gamma)].$$

Не является сложной задачей переход от (3) к уравнению с центрированной пуассоновской мерой $\tilde{v}(dt;d\gamma)$. В новом уравнении появятся дополнительные слагаемые при $dt$:

$$\int [(F(t;x(t) + g(t;\gamma)) - F(t;x(t))]\prod(d\gamma) + \int G(t;x(t) + g(t;\gamma);\gamma)\prod(d\gamma).$$

Как видим, различия в представлениях обобщенной формулы Ито-Вентцеля носят «технический» характер. Важней ограничения на коэффициенты различных представлений исходных уравнений, допускающие существование и единственность решений, осуществимость приведенных переходов.